\documentclass[times]{oupau-ppn}

\usepackage{amssymb,setspace}
\usepackage{hyperref}
\usepackage{caption,graphicx}
\usepackage{natbib}

\newtheorem{thm}{Theorem}

\newcommand{\R}{\mathbb{R}}
\newcommand{\Z}{\mathbb{Z}}

\newcommand{\nrm}[1]{\|#1\|}
\renewcommand{\(}{\left(}
\renewcommand{\)}{\right)}
\newcommand{\T}{\mathsf{T}}
\renewcommand{\L}{\mathsf{L}}
\newcommand{\Lsp}{\mathrm{L}}

\newcommand{\nx}{\nabla_x}
\newcommand{\nalpha}{\partial_\alpha}
\newcommand{\Null}{\mathcal{N}}

\newcommand{\eps}{\varepsilon}

\newcommand{\A}{\mathsf{A}}

\renewcommand{\H}{\mathsf{H}}
\newcommand{\D}{\mathsf{D}}

\newcommand{\be}[1]{\begin{equation}\label{#1}}
\newcommand{\ee}{\end{equation}}
\newcommand{\email}[1]{\emph{E-mail:} {\rm\textsf{#1}}}

\usepackage{color}   

\def\fref#1{(\ref{#1})}

\begin{document}
\runningheads{J. Dolbeault, A. Klar, C. Mouhot, and C. Schmeiser}{Hypocoercivity and a Fokker-Planck equation for fiber lay-down}

\title{Exponential rate of convergence to equilibrium for a model describing fiber lay-down processes}

\author{Jean Dolbeault\affil{a}\corrauth, Axel Klar\affil{b}, Cl\'ement Mouhot\affil{c}, Christian Schmeiser\affil{d}}

\address{\affil{a}Ceremade (UMR CNRS 7534), Universit\'e Paris-Dauphine, Place de Lattre de Tassigny, F-75775 Paris C\'edex 16, France.\\
\email{dolbeaul@ceremade.dauphine.fr},\\
\affil{b}Technische Universit\"at Kaiserslautern, Fachbereich Mathematik, E.~Schr\"o\-dinger Stra\ss e, D-67663 Kaiserslautern, Germany.\\
\email{klar@itwm.fhg.de},\\
\affil{c}University of Cambridge, DAMTP, Centre for Mathematical Sciences, Wilberforce Road, Cambridge CB3 0WA, UK.\\
\email{Clement.Mouhot@ens.fr},\\
\affil{d}Fakult\"at f\"ur Mathematik, Universit\"at Wien, Nordbergstra{\ss}e 15, 1090 Wien, Austria.
\email{Christian.Schmeiser@univie.ac.at}}

\corraddr{Jean Dolbeault: Ceremade (UMR CNRS 7534), Universit\'e Paris-Dauphine, Place de Lattre de Tassigny, F-75775 Paris C\'edex 16, France. \email{dolbeaul@ceremade.dauphine.fr}}

\received{\today}

\begin{abstract}
This paper is devoted to the adaptation of the method developed in \cite{DMS,Dolbeault2009511} to a Fokker-Planck equation for fiber lay-down which has been studied in \cite{MR2375288,MR2350004}. Exponential convergence towards a unique stationary state is proved in a norm which is equivalent to a weighted $\Lsp^2$ norm. The method is based on a micro / macro decomposition which is well adapted to the diffusion limit regime.
\end{abstract}


\maketitle

\section{Introduction}

The understanding of the shapes generated by the lay-down of flexible fibers onto a conveyor belt is of great interest in the production process of nonwovens textiles that find their applications, e.g., in composite materials like filters, textile and hygiene industry. In \cite{MR2350004} a stochastic model for the fiber lay-down process, \emph{i.e.}~for the generation of a fiber web on a conveyor belt has been presented. Taking into account the fiber motion under the influence of turbulence, the process can be described by a system of stochastic differential equations. We shall focus on a very simple situation, which does not take into account the movement of the belt, and gives rise to a Fokker-Planck equation. Some numerical results will also be considered at the end of this paper.

An important criterion for the quality of the web and the resulting nonwoven material is how the solution converges to equilibrium. In particular, the speed of convergence to the stationary solution is important. The faster this convergence is, the more uniform the produced textile will be. From a technological point of view, process parameters should be adjusted such that the speed of convergence to equilibrium is optimal.

The trend to equilibrium for solutions of kinetic equations has been investigated in many papers using entropy methods; see for example \cite{MR1787105,MR2562709}. A simplified approach has been suggested in \cite{DMS,Dolbeault2009511}. This trend to equilibrium for the Fokker-Planck equation for fiber lay-down under consideration in this paper has already been investigated using Dirichlet forms and operator semi-group techniques in \cite{MR2443261}; an ergodic theorem and explicit rates of convergence have been established. In the present paper we prove the convergence at an exponential rate towards a unique stationary state in a weighted $\Lsp^2$ norm by adapting the method developed in \cite{DMS,Dolbeault2009511} to the setting of non-moving belts.

\section{The model and main results}

In the melt-spinning process of nonwoven textiles, hundreds of individual endless fibers obtained by the continuous extrusion through nozzles of a melted polymer are stretched and entangled by highly turbulent air flows to finally form a web on a conveyor belt (see \cite{MR2350004} for more details). We describe the motion of an individual fiber, neglecting interactions with the others.

An arclength parametrization of the laid down fiber in a coordinate system following the conveyor belt is given by 
$x_0(t)\in\R^2$, $t\ge 0$. The tangent vector is denoted by $dx_0(t)/dt = \tau(\alpha(t))$ with $\tau(\alpha) = (\cos\alpha,\sin\alpha)$, $\alpha\in S^1 = \R/2\pi\Z$. Since the lay-down process is assumed to happen at the constant normalized speed 1 (equal to the spinning speed), $x_0(t)$ can also be interpreted as the position of the lay-down point at time $t$. If the conveyor belt moves with velocity $\kappa\,e_1$, the history of the lay-down point in the laboratory frame (as opposed to the conveyor belt frame) is given by $x(t) = x_0(t) + t\,\kappa\,e_1$, \emph{i.e.,}
\be{x}
\frac{dx}{dt} = \tau(\alpha) + \kappa\,e_1 \,.
\ee
It is a natural restriction that the speed of the conveyor belt cannot exceed the lay-down speed: $0\le \kappa\le 1$, since otherwise a stationary lay-down point would be impossible. The lay-down process can now be determined by prescribing the dynamics of the angle $\alpha(t)$, decribed as a stochastic process. It is driven by a deterministic force trying to move the lay-down point towards the equilibrium position $x=0$ and by a Brownian motion modeling the effect of the turbulent air flow:
\be{alpha}
d\alpha = -\,\tau^\bot(\alpha) \cdot \nabla V(x)\,dt + A\,dW\,,
\ee
where $W$ denotes a one-dimensional Wiener process, $A>0$ measures its strength relative to the deterministic forcing, $\tau^\bot = d\tau/d\alpha = (-\sin\alpha,\cos\alpha)$, and $V(x)$ is a potential such that $e^{-V}$ is integrable with the normalisation $\int_{\R^2} e^{-V}dx = 1$ and $\nabla V(0) = 0$.

The system \fref{x}--\fref{alpha} defines a stochastic process on $\R^2\times S^1$. The corresponding probability density $f(t,x,\alpha)$ satisfies the Fokker-Planck equation
\[
\partial_t f + (\tau + \kappa\,e_1)\cdot\nabla_x f - \partial_\alpha ( \tau^\bot\cdot\nabla_x V f + D\,\partial_\alpha f) = 0\,,
\]
(with the diffusivity $D = A^2/2$) which will be the object of our study. The analysis of the long time behaviour is considerably simplified in the case of a nonmoving conveyor belt:
\[
\kappa=0\,.
\]
We shall assume that this assumption holds true from now on. 
In this case the Fokker-Planck equation is written as an abstract ODE
\be{FP}
\partial_t f + \T f = D\,\L f\,,
\ee
with $\T f = \tau\cdot\nabla_x f - \partial_\alpha ( \tau^\bot\cdot\nabla_x V f)$ and $\L f = \partial_\alpha^2 f$.

It is easily seen that $F(x,\alpha) = e^{-V(x)}$ is an equilibrium solution of \fref{FP}, lying in the intersection of the null spaces of $\T$ and $\L$: $\T F = \L F = 0$.

A convenient functional analytic setting is introduced by the scalar product
\[
\langle f, g\rangle := \int_{\R^2\times S^1} fg\,d\mu\,,\quad d\mu(x,\alpha) := \frac{dx\,d\nu(\alpha)}{F(x,\alpha)}\,,\quad d\nu(\alpha) := \frac{d\alpha}{2\pi}\,,
\]
and by the associated norm $\|f\|^2=\langle f, f\rangle$. On the space $L^2(\R^2\times S^1, d\mu)$, the operator $\T$ is skew symmetric, and the operator $\L$ is symmetric and negative semi-definite. Thus, we have
\be{entropy}
\frac{d}{dt} \frac{\|f-F\|^2}{2} = D\,\langle \L f, f\rangle = -D\,\|\partial_\alpha f\|^2 \,.
\ee
This identity reveals the main difficulty in proving convergence to equilibrium. The decay to equilibrium seems to stop, as soon as $f$ is in the null space of $\L$ consisting of all $\alpha$-independent distributions. On the other hand, the decay equation \fref{entropy} does not make use of the action of the operator $\T$ and, in particular, of the fact that the equilibrium $F$ is the unique probability density in $\Null(\T) \cap \Null(\L)$. Any $\alpha$-independent distribution function $f$ is indeed unstable under the action of $\T$, unless $f=F$. Hence, convergence to the equilibrium can be expected and will be proven to be exponential. This is a so-called \emph{hypocoercivity} result as defined in \cite{MR2562709}. A recently developed approach \cite{DMS,Dolbeault2009511} for proving hypocoercivity in the abstract setting \fref{FP} will be applied with a special emphasis on the behaviour of the decay rate as $D\to 0$ and $D\to\infty$. It requires assumptions on the potential $V$, which have already been used in~\cite{DMS}:
\begin{itemize}
\item[(H1)] {\sl Regularity:\/} $V\in W^{2,\infty}_{\rm loc}(\R^2)$.
\item[(H2)] {\sl Normalization:\/} $\int_{\R^2}e^{-V}dx=1$.
\item[(H3)] {\sl Spectral gap condition:\/} there exists a positive constant $\Lambda$ such that
\[
\int_{\R^2}|\nx u|^2\,e^{-V}\,dx \ge \Lambda\int_{\R^2} u^2\,e^{-V}dx
\]
for any $u\in H^1(e^{-V}dx)$ such that $\int_{\R^2}u\,e^{-V}dx=0$.
\item[(H4)] {\sl Pointwise condition:\/} there exists $c_1>0$ such that \\
$|\nabla_x^2 V(x)|\leq c_1\,\big(1+|\nx V(x)|\big)$ for any $x\in\R^2$.
\end{itemize}
Roughly speaking, (H2) and (H3) require a sufficiently strong growth of $V(x)$ as $|x|\to\infty$, whereas (H4) puts a limitation on the growth behavior. This leaves room, however, for a large class of confining potentials including $V(x) = (1+|x|^2)^\beta$, $\beta\ge 1/2$. 

In \cite{DMS}, the additional pointwise condition 
$\Delta_x V(x)\leq \frac\theta2\,|\nx V(x)|^2+c_0$ with $c_0>0$ and $\theta\in (0,1)$ has been required. This is, however,
a consequence of (H4) by $\Delta_x V(x) \le\sqrt2\,|\nx^2 V(x)|$ and the Young inequality 
$\sqrt2\,c_1\,|\nx V| \le |\nx V|^2/4 + 2\,c_1^2$, with $\theta=1/2$ and $c_0 = c_1 + 2\,c_1^2$.

\begin{thm}\label{thm:decay} Let $f_0\in L^2(\R^2\times S^1,d\mu)$ and let (H1)--(H4) hold. Then, for every $\eta>0$, the solution of \fref{FP} subject to the initial condition $f(t=0) = f_0$ satisfies
\[
\|f(t)-F\| \le (1+ \eta)\,\|f_0-F\|\,e^{-\lambda t} \quad\mbox{with}\quad \lambda = \frac{\eta}{1+\eta} \frac{C_1\,D}{1+ C_2\,D^2}\,,
\]
where $C_1$ and $C_2$ are two positive constants which depend only on the potential $V$. \end{thm}
As a consequence, we have
\[
\lambda = O(D) \quad\mbox{as } D\to 0 \quad\mbox{and}\quad \lambda = O(D^{-1}) \quad\mbox{as } D\to\infty \,.
\]
Both results are sharp, as can be seen from the toy problem in \cite{DMS}. As $D\to 0$, dissipation is provided by the $O(D)$ right hand side of \fref{FP}, which dominates $\lambda$. On the other hand, the dynamics for large $D$ can be described by a macroscopic limit: see \cite{MR2375288}. When $D\to\infty$, the correct time scale is $t=O(D)$ and corresponds to a parabolic scaling, and therefore $\lambda = O(D^{-1})$ had to be expected.

However our method is not sharp in the sense that we cannot expect to obtain the optimal coefficients in the limiting cases above. Again this can be seen from the toy problem in \cite{DMS}, where the spectral gap can be explicitly computed.

\section{Proof of Theorem \ref{thm:decay}}

\subsection{The modified entropy}

We introduce the deviation $g := f-F$, satisfying \fref{FP} subject to $g(t=0) = f_0-F$.
Following \cite{DMS}, we denote the orthogonal projection to $\Null(\L)$ by
\[
\Pi g := \rho_g = \int_{S^1} g\,d\nu \,.
\]
Note that $\int_{\R^2} \rho_g\,dx = 0$ holds. In the following, we shall also need
\be{TPi}
\T\Pi g = \tau\cdot (\nx \rho_g + \rho_g \nx V) = \tau\cdot e^{-V}\nx(e^V \rho_g)\,,
\ee
with the consequence
\be{PiTPi}
\Pi\T\Pi = 0\,,
\ee
which is essential for the applicability of the method of \cite{DMS}. It implies that the macroscopic limit (corresponding to $D\to\infty$) in~\fref{FP} is diffusive: see \cite{MR2375288}.

With the help of the operator
\[
\A = \big(1+(\T\Pi)^*\T\Pi\big)^{-1}(\T\Pi)^*
\]
and an appropriately chosen $\eps>0$, the modified entropy functional is defined by 
\[
\H[g] := \frac{1}{2}\,\|g\|^2 + \eps\,\langle \A g,g\rangle\,.
\]
By \cite[Lemma 1]{DMS}, \fref{PiTPi} implies
\[
\|\A g\| \le \frac{1}{2}\,\|(1-\Pi)g\| \,.
\]
Hence the modified entropy functional is bounded from above and below by the square of the norm for any $\eps\in(0,1)$. More precisely we have:
\[
\frac{1-\eps}{2}\,\|g\|^2 \le \H[g] \le \frac{1+\eps}{2}\,\|g\|^2 \,.
\]
A straightforward computation gives
\be{dotH}
\frac{d}{dt} \H[g] = -\D[g]\,,
\ee
with the entropy dissipation functional
\be{diss}
\D[g] = -D\,\langle \L g, g\rangle + \eps\, \langle \A\T\Pi g,g\rangle + \eps\, \langle \A\T(1-\Pi)g,g\rangle - \eps\, \langle \T\A g,g\rangle - \eps\, D\,\langle \A\L g,g\rangle \,.
\ee

\subsection{Microscopic and macroscopic coercivity}
The first term on the right hand side of \fref{diss} has already been computed in \fref{entropy}. With $d\nu=d\alpha/(2\pi)$, the Poincar\'e inequality on $S^1$,
\[
\int_{S^1} |\nalpha g|^2\,d\nu \ge \int_{S^1} \(g - \textstyle{\int_{S^1} g\,d\nu} \)^2\,d\nu
\]
implies the \emph{microscopic coercivity} property
\be{coerc1}
-\langle \L g,g\rangle \ge \|(1-\Pi)g\|^2 \,.
\ee
The operator $\A\T\Pi = (1+(\T\Pi)^* \T\Pi)^{-1} (\T\Pi)^* \T\Pi$ shares its spectral decomposition with $(\T\Pi)^* \T\Pi$. For the latter we have, using \fref{TPi}
\[
\langle (\T\Pi)^* \T\Pi g,g\rangle = \|\T\Pi g\|^2 = \frac{1}{2}\,\int_{\R^2\times S^1} e^{-V} \left| \nx u_g\right|^2\,dx\,d\nu\,,
\]
with $u_g = e^V \rho_g$. The spectral gap condition (H3) implies the \emph{macroscopic coercivity} property
\[
\langle (\T\Pi)^* \T\Pi g,g\rangle \ge \frac{\Lambda}{2}\,\|\rho_g\|^2\,,
\]
leading to
\be{coerc2}
\langle \A\T\Pi g,g\rangle \ge \frac{\Lambda}{2+\Lambda}\,\|\Pi g\|^2 \,.
\ee
By \fref{coerc1} and \fref{coerc2}, the sum of the first two terms in the entropy dissipation \fref{diss} is coercive. This will also be sufficient for controlling the remaining three terms, if the operators $\A\T$, $\T\A$, and $\A\L$ are bounded, for $\eps>0$, small enough.

\subsection{Boundedness of auxiliary operators}
By \cite[Lemma 1]{DMS}, we know that
\be{coerc3}
\|\T\A g\| \le \|(1-\Pi)g\| \,.
\ee
The computation
\[
(\T\Pi)^* \L g = -\Pi\T \L g = -\nx\cdot \Pi(\tau\,\nalpha^2 g) = \nx\cdot\Pi(\tau\,g) = -(\T\Pi)^* g
\]
shows that $\A\L = -\A$ and, thus,
\be{coerc4}
\|\A\L g\| \le \frac{1}{2}\,\|(1-\Pi)g\| \,.
\ee
The most elaborate part of the analysis is to prove the boundedness of $\A\T$. Following the approach of \cite{DMS}, we consider its adjoint
\[
(\A\T)^* = -\T^2\Pi \big(1+(\T\Pi)^*\T\Pi\big)^{-1}\,.
\]
For a given $g\in L^2(\R^2\times S^1,d\mu)$, we introduce $h = \big(1+(\T\Pi)^*\T\Pi\big)^{-1} g$ which, after solving for $g$ and applying~$\Pi$, becomes
\be{elliptic}
\rho_g = \rho_h - \Pi \T^2 \rho_h = e^{-V}u_h - \frac{1}{2}\,\nx\cdot (e^{-V} \nx u_h)
\ee
with $u_h = e^V \rho_h$. A straightforward computation gives
\[
(\A\T)^* g = -\T^2 \rho_h = e^{-V} \left[ (\tau\cdot\nx)^2\,u_h - (\tau^\bot\cdot \nx V)\,(\tau^\bot \cdot \nx u_h) \right]
\]
and, as a consequence,
\[
\|(\A\T)^* g\| \le \left\|\nx^{\otimes 2} u_h\right\|_{L^2(\R^2,e^{-V}dx)} + \left\| |\nx V|\, |\nx u_h|\right\|_{L^2(\R^2,e^{-V}dx)}\,.
\]
Therefore, in order to prove the boundedness of $(\A\T)^*$ (and, thus, of $\A\T$), we need to prove the boundedness of the right hand side in terms of $\|\rho_g\|$ for the solution $u_h$ of the elliptic equation \fref{elliptic}. This $L^2\to H^2$ (with weight~$e^{-V}$) elliptic regularity result has been derived in \cite{DMS} under the assumptions (H1)--(H4) (see Proposition~5 for the first term and Lemma~8 for the second). Collecting these results gives $\|(\A\T)^* g\| \le C_V\,\|g\|$ and therefore
\be{coerc5}
\|\A\T(1-\Pi)g\| \le C_V\,\|(1-\Pi)g\|\,,
\ee
where $C_V$ depends only on the potential.

\subsection{Hypocoercivity}
Inserting estimates \fref{coerc1}--\fref{coerc4} and \fref{coerc5} in \fref{diss} gives
\begin{eqnarray}
\D[g] &\ge& D\,\|(1-\Pi)g\|^2 + \tfrac{\eps\,\Lambda}{2+\Lambda}\,\|\Pi g\|^2 - \eps \(C_V+1+\tfrac{D}{2}\)\|(1-\Pi)g\|\,\|g\| \nonumber\\
&\ge& \(D - \eps \(C_V+1+\tfrac{D}{2}\)\)\|(1-\Pi)g\|^2 + \tfrac{\eps\,\Lambda}{2+\Lambda}\,\|\Pi g\|^2-\, \eps \(C_V+1+\tfrac{D}{2}\)\|(1-\Pi)g\|\,\|\Pi g\| \label{diss-est}\\
&\ge& \(D - \eps \(C_V+1+\tfrac{D}{2}\)\(1+\tfrac{1}{2\delta}\)\)\|(1-\Pi)g\|^2 +\, \eps\( \tfrac{\Lambda}{2+\Lambda} - \tfrac{\delta}{2}\(C_V+1+\tfrac{D}{2}\)\) \|\Pi g\|^2\,, \nonumber
\end{eqnarray}
for an arbitrary $\delta>0$. This shows already that coercivity can be achieved by first choosing $\delta$ and then $\eps$, both small enough. With the choice
\[
\delta = \frac{\Lambda}{(2+\Lambda)\,(C_V+1+D/2)}\,,
\]
the coefficients on the right hand side of \fref{diss-est} can be written as $D-\eps\, \mathsf r(D)$ and $\eps\,\mathsf s$ with
\[
\mathsf r(D) := \frac{1}{2\,\Lambda} \(2\Lambda + (2+\Lambda)\(C_V+1+\tfrac{D}{2}\)\)\(C_V+1+\tfrac{D}{2}\)\,,\quad \mathsf s := \frac{\Lambda}{2\,(2+\Lambda)}\,.
\]
Then the optimal choice of $\eps$, considering the form of the coefficients, would be
\[
\overline\eps(D) := \frac{D}{\mathsf r(D)+s} \,.
\]
However, we also have to guarantee $\eps<1$ for the definiteness of $\H[g]$ and actually, even stronger, $\frac{1+\eps}{1-\eps} \le (1+\eta)^2$ will be needed below, which can be guaranteed by the requirement $\eps \le \frac{\eta}{1+\eta}$. Moreover the two conditions are equivalent at first order for $\eta>0$, small. These considerations lead to the choice
\[
\eps = \frac{\eta}{1+\eta}\,\frac{\overline\eps(D)}{\overline\eps_{max}}\,,\quad\mbox{with } \overline\eps_{max} = \max\big\{ 1, \max_{D>0} \overline\eps(D)\big\}\,,
\]
which is finite because of $\overline\eps(0) = \overline\eps(\infty) = 0$. With this choice,
\[
D-\eps\, \mathsf r(D) \ge \eps\, \mathsf s \ge \frac{\eta}{1+\eta}\,\frac{2\,C_1\,D}{1+ C_2\,D^2} =: 2\,\lambda\,,
\]
with appropriately chosen constants $C_1$, $C_2>0$, depending only on $\Lambda$ and $C_V$ and, thus, only on the potential~$V$. The estimate
\[
\D[g]\ge 2\,\lambda\,\|g\|^2 \ge \frac{4\,\lambda}{1+\eps}\,\H[g] > 2\,\lambda\,\H[g]
\]
follows. Using this in \fref{dotH} and the Gronwall lemma imply
\[
\H[f(t)-F] \le \H[f_0-F]\,e^{-2\lambda t} \,.
\]
Finally we obtain for the norm
\begin{eqnarray*}
\|f(t)-F\|^2 &\le& \frac{2}{1-\eps}\, H[f(t)-F] \le \frac{2}{1-\eps}\, H[f_0-F]\,e^{-2\lambda t} \\
&\le& \frac{1+\eps}{1-\eps}\, \|f_0-F\|^2\,e^{-2\lambda t} \le (1+\eta)^2\,\|f_0-F\|^2\,e^{-2\lambda t}\,,
\end{eqnarray*}
which completes the proof of Theorem~\ref{thm:decay}.


\section{Concluding remarks}

\subsection{Numerical investigations}

It is interesting to compare the rates predicted by the above results, which are only upper bounds, with numerical rates of convergence. We use a classical Monte-Carlo method with an Euler-Maruyama discretization scheme for all computations. A numerical investigation of the equations using a semi-Lagrangian method can be found in \cite{MR2475656}.

\begin{figure}\begin{center}
\includegraphics[clip,width=0.75\columnwidth]{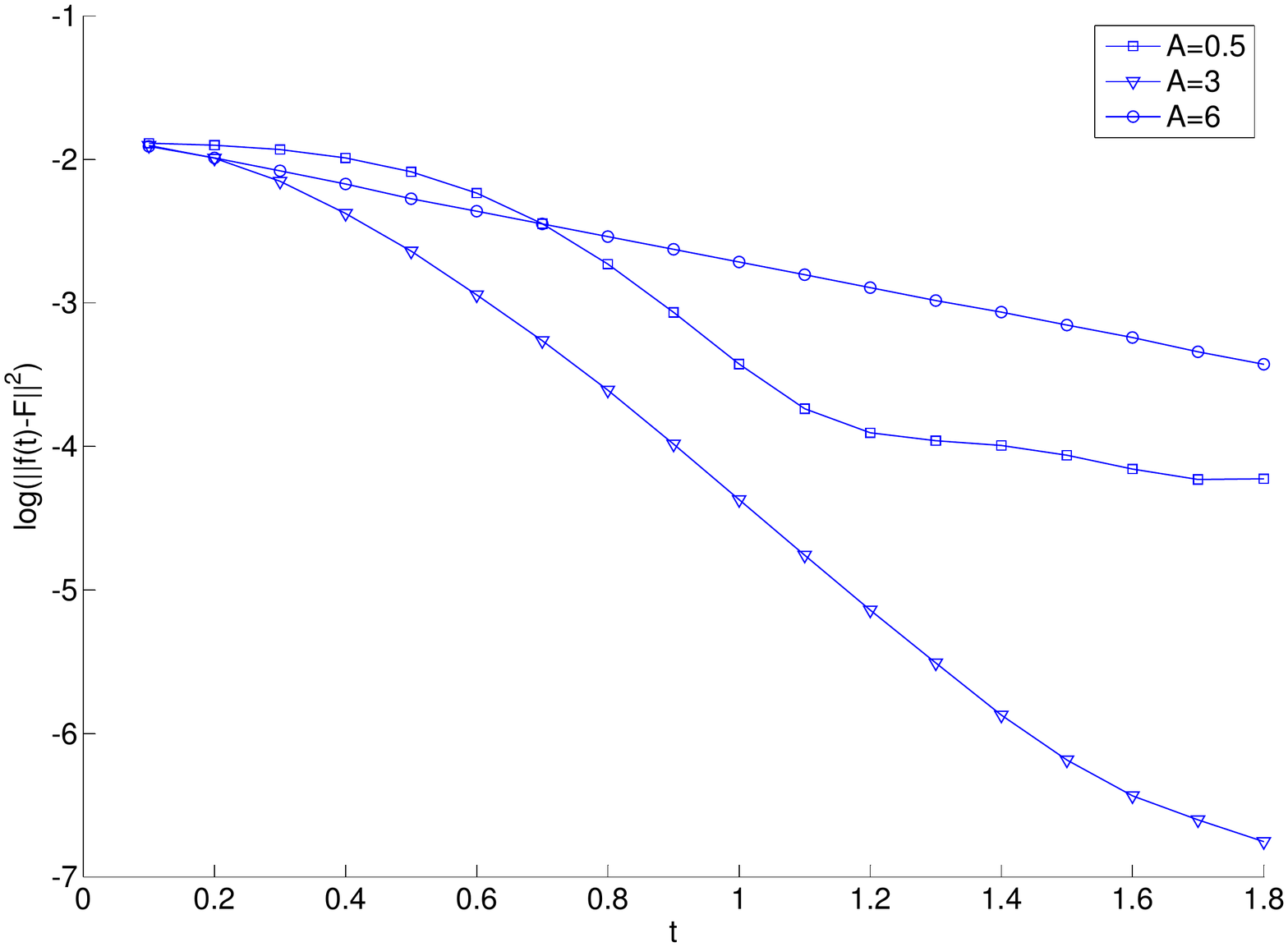}
\end{center}
\caption{\textit{Plot of $t\mapsto\log\nrm{f(t,\cdot,\cdot)-F}^2$ for $A = 0.5,3,6$. }}
\label{figentropy}\end{figure}

\begin{figure}\begin{center}
    \includegraphics[angle=90,clip,width=0.75\columnwidth]{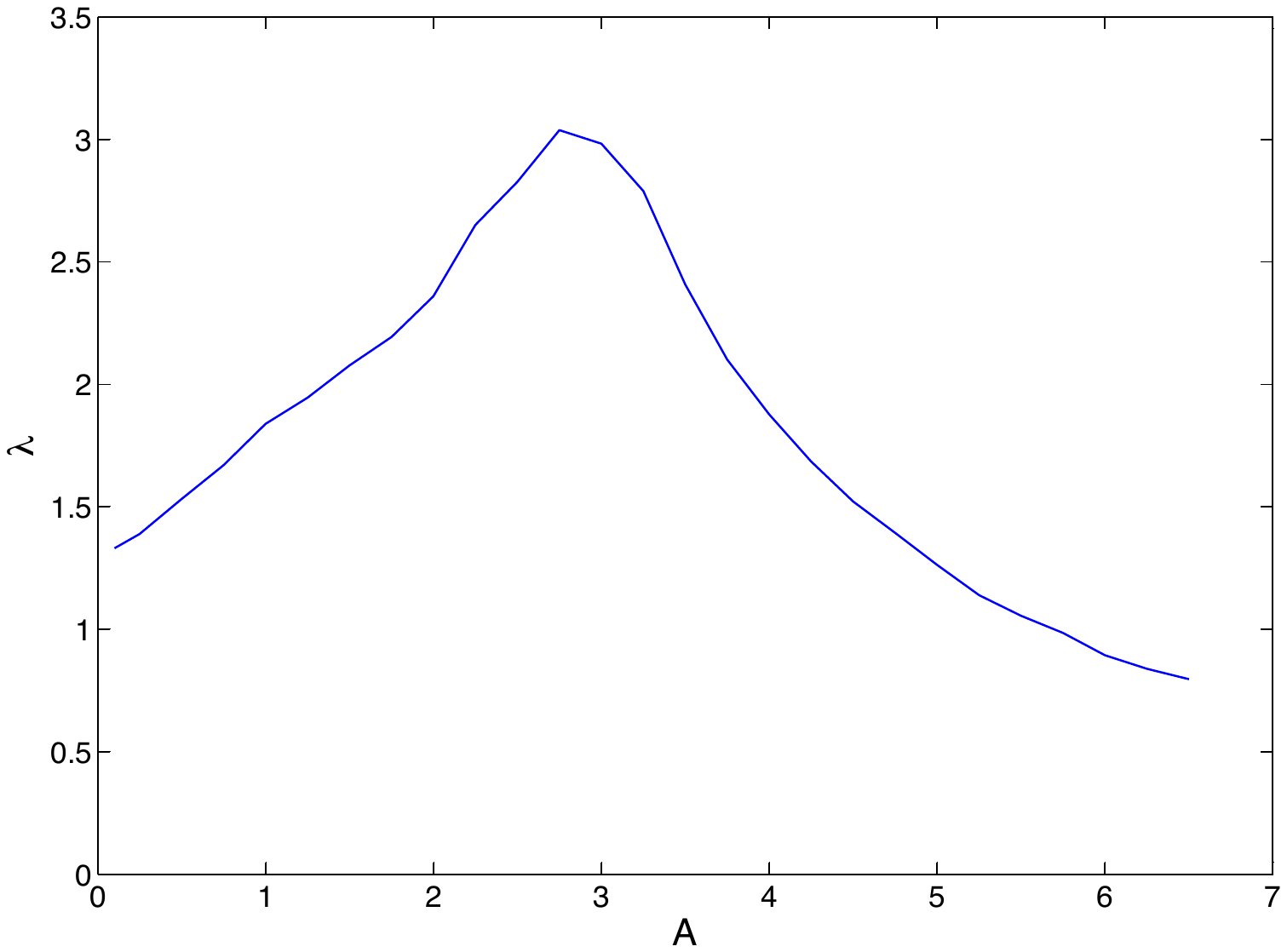}
\end{center}
\caption{\textit{Plot of $\lambda = \lambda (A)$ for different values of $A$.}}
\label{figlambda}\end{figure}

The exponential decay of the $L^2$-difference to the stationary solution is observed in Figure \ref{figentropy}. In Figure \ref{figlambda} the decay rates $\lambda$ have been obtained from the above simulations for various values of $A$ using a least square fit. The rate given by Theorem~\ref{thm:decay}, \emph{i.e.} $\lambda \sim \frac{C_1\,D}{1\,+\, C_2\,D^2}$, $D=A^2/2$, fits qualitatively very well the curve obtained in Figure~\ref{figlambda} when $A$ is away from $0$. In particular, values of $A$ with an optimal rate of convergence can be determined from the numerical as well as the analytical results.

\subsection{Perspectives}

\begin{enumerate}

\item \emph{Models where stationary solutions are not known explicitly}. The model considered in this paper can be extended in different directions, for instance by taking into account the movement of the belt, or by models where fibers have smoother trajectories than the ones considered in this paper, see \cite{MR2525723,MR2375288}. In these cases the stationary solutions are not always known explicitly. The application of the entropy method presented above is then an open problem.

\item \emph{3-D models}. To model fluid flow through a fiber web, the model has to be extended to three dimensional situations, see \cite{KMW09}. For such a model exponential convergence to equilibrium, at least for the case $\kappa=0$, can be proven with the same methods as in this paper.

\end{enumerate}

\medskip\begin{minipage}{163.5mm}\linespread{0.5}\selectfont\small\noindent{\sl Acknowledgements.\/} {\small This research project has been supported by the ANR project \emph{CBDif-Fr} and \emph{EVOL}, by the Excellence Center for Mathematical and Computational Modeling (CM)$^2$ and by Deutsche Forschungsgemeinschaft (DFG), KL 1105/18-1. The authors thank K. Fellner and P. Markowich for the organization of a conference on \emph{Modern Topics in Nonlinear Kinetic Equations} (DAMTP, Cambridge, April 20-22, 2009) were this research project was initiated.}\end{minipage}\bigskip

\bibliographystyle{plainnat}
\bibliography{References}
\end{document}